\documentclass[xypic,amscd,syntonly,amssymb,verbatim,12pt]{amsart}

\usepackage{graphicx}
\usepackage[all]{xy}
\usepackage{amssymb}
\usepackage{amsmath}
\usepackage[mathscr]{euscript}

\usepackage{epsfig}
\theoremstyle{plain}
\newtheorem{Thm}{Theorem}
\newtheorem{Prop}[Thm]{Proposition}
\newtheorem{Cor}[Thm]{Corollary}

 \theoremstyle{definition}
\theoremstyle{remark}

\numberwithin{equation}{section}

\begin{document}
  \title{Charge of $D$-branes on singular varieties}

 \author{ ANDR\'{E}S   VI\~{N}A}
\address{Departamento de F\'{i}sica. Universidad de Oviedo.    Garc\'{\i}a Lorca 18.
     33007 Oviedo. Spain. }
\email{vina@uniovi.es}

 \maketitle
\begin{abstract}

Considering  the $D$-branes on a variety $Z$ as the objects of the derived category $D^b(Z)$,  we propose a definition for the charge of $D$-branes on not necessarily smooth varieties. 
We define the charge 
$Q({\mathcal G})$  of 
 ${\mathcal G}\in D^b(Z)$ 
 as an element of the homology of $Z$, so that the mapping $Q$
is compatible with the pushforward by proper maps between varieties.

Given a generic anticanonical hypersurface $Y$ of a toric variety $X$ defined by a reflexive polytope, we 
 express the charge of a line bundle on $Y$ defined by a divisor $D$ of $X$
 in terms of  intersections of $D$ with cycles  determined by the   polytope faces.


\end{abstract}

\section {Introduction} \label{S:intro}

 The $D$-branes of type $B$ on a Calabi-Yau variety $Z$ are the objects of $D^b(Z)$, the   bounded derived category of algebraic coherent sheaves on $Z$ (see monograph \cite{Aspin-et}). Given a brane ${\mathcal G}$ on $Z$, when $Z$ is a {\em smooth} variety, the charge of ${\mathcal G}$ 
is the element of the cohomology $H^*(Z)$  defined 
by the cup product of the Chern character of ${\mathcal G}$ and the root square of the Todd class of $Z$ \cite{Aspin}
\begin{equation}\label{charge:0}
{\rm ch}({\mathcal G}) \sqrt{{\rm td}(Z)}.
 \end{equation}
For some branes ${\mathcal G}$,  the integration of this cohomology class gives the index of a differential  Dirac operator.

Obviously,  when the variety $Z$ is singular its Todd class is not defined. On the other hand, in the singular case there are not locally free resolutions of the coherent sheaves, so 
the  Chern character for vector bundles
 has not a direct extension to coherent sheaves.  
 
 In Section \ref{S:Charge} of this note, we propose a definition for the charge of a brane on a  projective variety $Z$. This charge will be an element of $H_*Z$,  the rational homology of $Z$. When $Z$ is not necessarily smooth, it is not possible to demand  the charge to be related with the index of some differential operator. However,  
the simple situation that we next describe  suggests a property to be satisfied by a charge definition.

Given a closed  subvariety $Y$ of   $X$, 
let $i$ denote   the closed embedding $i:Y\hookrightarrow X$. The extension by zero of a coherent sheaf ${\mathcal F}$ on $Y$ is the coherent sheaf $i_{!}{\mathcal F}$ on $X$.  With respect to a definition of charge, it seems natural to demand the ``equality'' between the charge of ${\mathcal F}$ and  the charge of its extension. The  ``equality'' must be understood in the following sense  
\begin{equation}\label{i_*Q}
i_*Q({\mathcal F})=Q(i_!{\mathcal F}),
 \end{equation} 
where $Q(\,\star\,)$ denotes the charge of the brane  $\star$ and $i_*:H_*Y\to H_*X$ is the homomorphism induced by the inclusion $i$. We call the equality (\ref{i_*Q}) the consistency of the charge with respect to the pushforward.

The approach consisting of embedding a singular space in a smooth one and extending the objects to the smooth space
has been  used  in \cite{B-F-M} to prove the Grothendieck-Riemann-Roch (G-R-R) theorem for singular varieties, and in \cite{Hartshorne} to define an algebraic de Rham cohomology for schemes.

 On the other hand, the above consistency condition for the charge  leads us to the G-R-R theorem \cite{B-S}.  So, we will translate to the language of branes the 
  process  carried out in the proof of the G-R-R theorem for singular varieties in
\cite{B-F-M}.
 The idea consists of embedding $Z$ in a smooth variety $M$ and    extending to $M$ the branes given on $Z$. 
The next step is to localize on $Z$ the Chern character of the extended branes, via the relative Grothendieck groups. 

 We propose  the localization of the Chern character of the extended branes, by means of gauge fields on  the smooth variety in which $Z$ is embedded. More precisely, let $j:Z\hookrightarrow M$ be an embedding of $Z$ in the smooth variety $M$ and ${\mathcal F}$ a coherent sheaf over $Z$.
 \begin{enumerate}
 \item Since $M$ is smooth, there is a resolution $E_{\bullet}$ of $j_!{\mathcal F}$ by locally free sheaves. 
 \item By the additivity of the Chern character, ${\rm ch}(j_!{\mathcal F})$ equals the character of the resolution.
 \item The restriction to $M\setminus Z$ of $E_{\bullet}$ is an exact sequence. So, 
 $$\sum_i (-1)^i{\rm ch}(E_i)$$
 restricted to $M\setminus Z$ is represented by an exact form.
 \item The Chern character of $E_i$ can be defined by the closed form 
	${\rm trace}({\rm exp}(k\Theta_i)),$ where $\Theta_i$ is the strength of any gauge field on $E_i$, and $k= \tfrac{\sqrt{-1}}{2\pi}$. 
	 \item By $(3),$  the restriction of the closed form
	$$\sum (-1)^i {\rm trace}({\rm exp}(k\Theta_i))$$
	to $M\setminus Z$  is an exact form. Hence, it defines an element of the de Rham relative cohomology $H^*(M,\,M\setminus Z)\simeq H_*(Z)$, that we consider as the localization of ${\rm ch}j_!{\mathcal F}$ on $Z$.
  This homology class is denoted in \cite{B-F-M}  by ${\rm ch}^M_Z({\mathcal F})$.
 \end{enumerate}
 We will define the  charge $Q({\mathcal F})$ of the brane ${\mathcal F}$ by the cap product of the Todd class of $M$ with   localized homology class ${\rm ch}^M_Z({\mathcal F})$.
 
 When $Y$ is a local complete intersection in a smooth variety and the brane  is a complex of locally free sheaves,  its  charge 
 can be expressed in terms of the character of the brane and the Todd class of the virtual tangent bundle of $Y$. The result is stated in Proposition \ref{Prop:Charge}.

 As a consequence of the consistency of $Q$ with  pushforwards, the arithmetic genus $\chi(Z,\,E)$ of a vector bundle $E$ over the variety $Z$
is the   $0$-degree component in the cap product ${\rm ch}(E)\cap Q({\mathcal O}_Z)$. This is in fact the statement of the generalized Hirzebruch-Riemann-Roch theorem \cite{B-F-M}.

Associated  to the variety $Z$ is the  exterior differential operator
$$d:\Omega^{\bullet}_{Z/{\mathbb C}}\rightarrow \Omega^{\bullet+1}_{Z/{\mathbb C}},$$
where $\Omega^{\bullet}_{Z/{\mathbb C}}$ is the algebraic de Rham complex, consisting of the sheaves of algebraic differentials on $Z$.

The algebraic de Rham cohomology $H^*_{dR}(Z{/{\mathbb C}})$ of $Z$ is by definition the hypercohomology ${\mathbb H}^*(Z,\, 
\Omega^{\bullet}_{Z/{\mathbb C}})$ \cite{Groth}
$$H^*_{dR}(Z {/{\mathbb C}})={\mathbb H}^*(Z,\, 
\Omega^{\bullet}_{Z/{\mathbb C}}).$$
 When $Z$ is {\it smooth}, one can consider $Z^{\rm an}$, the underlying complex manifold. Grothendieck proved that the singular cohomology of 
 $Z^{\rm an}$ is isomorphic to the algebraic de Rham cohomology
\begin{equation}\label{Groth-Iso}
H^*(Z^{\rm an},\,{\mathbb C})\simeq H^*_{dR}(Z {/{\mathbb C}}). 
\end{equation}
But, there are examples of {\it singular} varieties for which (\ref{Groth-Iso}) does not hold. Thus, in general, the algebraic de Rham cohomology groups have not a clear topological interpretation.

It may be tempting to think about a possible relation between $Q$ and the index of the operator $d$. But the following two observations discard this relation.
(a) According to the above mentioned property, the arithmetic genus $\chi(Z,\,{\mathcal O}_Z)$ of $Z$ is a topological invariant determined by $Q({\mathcal O}_Z)$. (b) In general, ${\rm dim}\,H^k_{dR}(Z {/{\mathbb C}}) $ is not a topological invariant of $Z$.


\smallskip 

 In Section \ref{S:Euler},  we consider 
  branes on  Calabi-Yau anticanonical hypersurfaces of the smooth toric variety determined by a reflexive polytope.

 Let $\Delta$ be a reflexive polytope in $({\mathbb R}^n)^*$ that defines the toric variety $X$, which we assume smooth.  Batyrev proved that there are anticanonical hypersurfaces of $X$ that are Calabi-Yau varieties (see \cite{Batyrev}). Let $Y$ denote such a hypersurface.
 Let $E$ be  the restriction to $Y$ of a line bundle over $X$ defined by a divisor $D$. By applying the above mentioned Hirzebruch-Riemann-Roch theorem, when ${\rm dim}_{\mathbb C}\, X=4$, we express the arithmetic genus $\chi(Y,\,E)$ of the line bundle $E$, in terms of intersections of homology classes of $X$ defined    by the faces   of the polytope $\Delta$ and by $D$ (see Theorem \ref{Thm:3}).

When ${\rm dim}_{\mathbb C}\, X=3$,  we   calculate the charge of a brane on $Y$ which is the restriction of a line bundle over $X$. The result is stated in Proposition \ref{Prop:charge}. In the case that $X$ is a complex surface, the charge of a brane on $Y$ defined by a divisor $D$ of $X$ can be expressed by a simple formula, in which only are involved $D$ and the canonical divisor of $X$ (Corollary \ref{Coro:charge}).

 The multiple $(n-1)(-K_X)$ of the anticanonical divisor of $X$ determines a line bundle $L'$ on $X$. We denote by $L$  for the restriction of $L'$   to $Y$. In \cite{Vina18}, we proved
that the family $\{{ L}^{\otimes k}\,|\,k\in{\mathbb Z}\}$ is a spanning class of the category $D^b(Y)$ \cite{Bridgeland, Huybrechts}.

In Corollary \ref{Cor:L}, we give an expression for $\chi(Y,\,L)$, when $X$ is a four dimensional variety. 
The value of the charge  $Q(L)$ of the brane $L$ is given in Corollary \ref{Cor:mathacalL}, 
assuming that  ${\rm dim}_{\mathbb C}\, X=3$.


\section{Charge of a brane}\label{S:Charge}

 \subsection{Grothendieck groups.} Firstly, we will review some properties of the Grothendieck groups (for more  details   see \cite{Fulton-Int,Manin}) and introduce some notations. Given a variety $Z$,  we denote by  $\mathfrak{Coh}(Z)$  the category of coherent sheaves on $Z$. The Grothendieck group $K_0Z$ is 
 the abelian group generated by the objects of $\mathfrak{Coh}(Z)$ modulo the relations
 ${\mathcal F}={\mathcal F'}+{\mathcal F''}$, when there is an exact short sequence
$$0\to{\mathcal F'}\to{\mathcal F}\to{\mathcal F''}\to 0$$
in the category $\mathfrak{Coh}(Z)$.

 Let $({\mathcal F}^{\bullet},\,d^{\bullet})$ be a bounded complex consisting of coherent sheaves over $Z$. From the exactness of the short sequences of $\mathfrak{Coh}(Z)$
 $$0\to{\mathcal Im}\,d^{i-1} \to{\mathcal Ker}\,d^{i}\to {\mathcal H}^i({\mathcal F}^{\bullet})\to 0 $$
 $$ 0\to {\mathcal Ker}\,d^{i}\to \mathcal{F}^i\to{\mathcal Im}\,d^{i}\to 0,$$
 it follows the following equality between elements of $K_0Z$
 $$\sum_k(-1)^k [{\mathcal F}^k] = \sum_k(-1)^k [{\mathcal H}^k({\mathcal F}^{\bullet})].$$
 Thus, 
 \begin{equation}\label{DmapsstoK}
 \alpha_Z:D^b(Z)\to K_0Z,\;\;\;\;\;{\mathcal G}\mapsto [{\mathcal G}]:= \sum_k(-1)^k [{\mathcal G}^k]
 \end{equation}
  is a well-defined map.

 By $K^0Z$ we denote the Grothendieck ring of algebraic vector bundles on $Z$, where the ``multiplication" is defined through the tensor product. As the Chern character ${\rm ch}$ is additive in exact sequences, it descends to a ring homomorphism from $K^0Z$ to the rational cohomology of $Z$, 
${\rm ch}:(K^0Z)\otimes {\mathbb Q}\to H^*Z$.

  The natural homomorphism $K^0Z\to K_0Z$, together with the fact that tensoring with locally free sheaves is an exact functor, defines a $K^0Z$-module structure on $K_0Z$. That is, given $[E]\in K^0Z$ and  $[{\mathcal F}]\in K_0Z$,
$$[E]\cdot[{\mathcal F}]=[E\otimes_{{\mathcal O}_Z} {\mathcal F}].$$

 The pullback of vector bundles by a map $g:Z\to W$, defines a homorphism of rings $g^*:K^0W \to K^0Z$. 
	Similarly, given   a proper map between the varieties $f:Z\to W$ and a coherent sheaf ${\mathcal H}$ over $Z$, the pushforward
	$${\mathcal H}\mapsto \sum_k (-1)^k[R^kf_*{\mathcal H}],$$
	determines a homorphism $f_*:K_0Z \to K_0W.$  Thus, $K_0$ can be considered as covariant functor defined on the category whose objects are  algebraic varieties, and whose arrows are proper morphisms. The morphisms $f^*$ and $f_*$ are related by the projection formula
	$$f_*(f^*(a )\cdot b)=a\cdot f_*(b),$$
	for $a\in K^0W$ and $b\in K_0Z$.

When $Z$ is {\em smooth}, any coherent sheaf on $Z$ admits a finite locally free resolution. From this property, it follows that $K^0Z$ and $K_0Z$ can be identified.
Thus,  in this case, the Chern character of a brane ${\mathcal G}$ is the image of ${\mathcal G}$ 
by the composition  
 \begin{equation}\label{SmoothCaracter}
D^b(Z)\overset{\alpha_Z}{\longrightarrow} K_0Z\simeq K^0 Z\overset{\rm ch}{\longrightarrow} H^*(Z;\,{\mathbb Q}).
\end{equation}


\smallskip
\subsection{Consistency of $Q$ with the pushforward.}\label{Ss:Consistence}
If, in the smooth case, we define the charge of a brane as its Chern character according to (\ref{SmoothCaracter}), this definition does not satisfy the consistency condition (\ref{i_*Q})   mentioned in the Introduction. 
The following example 
  suggests a   definition of  charge  for which  (\ref{i_*Q}) holds.  This example may be illustrative only for  readers, who are not familiar with the Grothendieck-Riemann-Roch theorem \cite{B-S}; nevertheless, we  deduce some results and introduce notations which will be used later.

 Let $i:Y\hookrightarrow X$  be a closed smooth subvariety of the  smooth compact variety $X$. We assume that $Y$ is the divisor defined by a section $s$ of a line bundle $N$ over $X$. Let $F'$ be an algebraic vector bundle over $X$ and set $F$ for the pullback 
  by $i$ of the $F'\otimes N$. We denote by ${\mathcal F}$ the sheaf over $Y$ determined by $F$. 
As $R^ki_*{\mathcal F}=0$ for $k>0$, $i_*{\mathcal F}$ reduces to $R^0i_*{\mathcal F}$, i.e. to the extension  $i_!{\mathcal F}$
  of the sheaf ${\mathcal F}$  to $X$ by zero. 
 
 One has the   morphism $F'\to F'\otimes N$ defined by
$$\sigma\in F'(V)\mapsto \sigma\otimes s|_V\in (F'\otimes N)(V)\,\;\hbox{for every open}\; V\subset X.$$
As $s$ is a global section of $N$ that does not vanishes identically over any open of $X$ and is never zero on the complement of Y, we have
  the following exact sequence of sheaves
 \begin{equation}\label{locallyfree}
 0\to F' \to F'\otimes N\to i_!{\mathcal F}\to 0,
  \end{equation}
which in turn defines  a locally free resolution of $i_!{\mathcal F}$.
Hence,
\begin{equation}\label{rmch}
{\rm ch}(i_!{\mathcal F})={\rm ch}(F'\otimes N)-{\rm ch}(F')= {\rm ch}(F')({\rm e}^a-1),
\end{equation}
where $a$ is the first Chern class of $N$.

On the other hand, one has the Gysin homomorphism between the cohomologies $\iota_*:H^*Y\to H^*X$ satisfying $\iota_* i^*\alpha =a\alpha$, for $\alpha\in H^*X$. Then, by definition of $F$
\begin{equation}\label{iota_*}
\iota_*{\rm ch}{ F}=\iota_*i^*({\rm ch}(F'\otimes N))=a\,{\rm e}^a{\rm ch}(F').
\end{equation}
Thus, $\iota_*{\rm ch}{ F}\ne {\rm ch}(i_!{\mathcal F})$. Since $i_*=P_X\circ\iota_*\circ P_Y^{-1}$, where $P_X$ and $P_Y$ are the corresponding Poincar\'e isomorphisms from cohomology to homology, it follows 
$$i_*(P_Y{\rm ch}{ F})\ne P_X{\rm ch}(i_!{\mathcal F}).$$
That is, the Poincar\'e dual of the Chern character does not satisfy condition (\ref{i_*Q}).

The G-R-R theorem \cite{B-S} shows us the way to overcome this difficulty.   From the exact sequence of bundles
$$0\to T_Y\to T_X|_Y\to N|_Y\to 0,$$
together with the fact that the Todd class is multiplicative on exact sequences, it follows
$${\rm td}(Y)=i^*\big(\frac{1-{\rm e}^{-a}}{a}\big)i^*({\rm td}(X)).$$ 
By (\ref{iota_*})
$$\iota_*({\rm ch}({ F})\,{\rm td}(Y))=\iota_*i^*\big({\rm ch}(F')\,{\rm e}^a\,\frac{1-{\rm e}^{-a}}{a}\, {\rm td}(X)   \big)=
{\rm ch}(F')({\rm e}^a-1)\,{\rm td}(X).$$
From (\ref{rmch}), it follows
$$\iota_*({\rm ch}({\mathcal F})\,{\rm td}(Y))={\rm ch}(i_!{\mathcal F})\,{\rm td}(X).$$

Therefore, if ${\mathcal G}$ is a coherent sheaf on a {\it smooth} compact variety $Z$,  it is reasonable to define its charge $Q({\mathcal G} )$ as the element of $ H_*(Z)$
\begin{equation}\label{Q({math G})}
Q({\mathcal G})=P_Z\big({\rm ch}({\mathcal G})\,{\rm td}(Z)\big).
\end{equation}

  It is possible, of course, to consider other definitions for the charge of a brane over a smooth variety. In particular, for complexes of vector bundles 
	on which are defined an elliptic differential operators, one can consider charges related with the index of those operators; this is the case of (\ref{charge:0}).


\smallskip
		\subsection{ Localization of Chern character.}
		As we said in Section \ref{S:intro}, to define the charge of a brane over a singular variety, we embed the variety in a smooth one and extend the brane given to a brane over the smooth variety. Then it is necessary to localize the Chern character of the extended brane. In the following paragraphs (a), (b) and (c) we show the ideas behind this localization.
	
	\smallskip
		
{\it (a) Locally complete intersection.}
Let  $i:Y\hookrightarrow X$ be a closed embedding of the subvariety $Y$ in the smooth variety $X$, and let us assume that $Y$ is a local complete intersection defined by a section $s$ of an algebraic  line vector bundle $N$ on $X$. That is, we are in the situation considered in Section \ref{Ss:Consistence}  {\em without assuming that $Y$ is smooth}. We will use the notations introduced in that section.

As $i_!{\mathcal F}$ vanishes out $Y$,  the resolution (\ref{locallyfree}) gives an isomorphism of sheaves between 
$F'$ and $F'\otimes N$, when we restrict them to $X\setminus Y$. As the section $s$ does not vanish on $X\setminus Y$, one can define a connection $\Hat\nabla$ on $N|_{X\setminus Y}$ by setting $\Hat\nabla s=0$.

 To extend this connection to $N$, we choose an  open neighbourhood $U$ of $Y$ (of course $U$ can be taken so small as one wishes). By means of a partition of unity subordinate to the covering $\{U,\, X\setminus Y\}$ of $X$ and an arbitrary connection on $N|_U$, we obtain a connection on $N$, that will be also denoted by $\Hat\nabla$. The curvature form $\Hat\Theta$ of $\Hat\nabla$  vanishes out $U$ and  the cohomology class defined by $\Hat\Theta$  satisfies 
 $$[\Hat\Theta] =\frac{2\pi}{\sqrt{-1}}\,c_1(N).$$
  That is, 
   $\frac{\sqrt{-1}}{2\pi}\,\Hat\Theta$ is a Thom form of the embedding $Y\hookrightarrow X$.
   
   It is well-known that if $E$ is a vector bundle of rank $l$ over a smooth variety $Z$ and $\nabla$ is a connection on $E$ with curvature $\Theta$, then 
		\begin{equation}\label{chern-Nabla}
		{\rm ch}(\nabla):=l+\sum_{k>0}\frac{1}{k!}\Big(\frac{\sqrt{-1}}{2\pi}\Big)^k{\rm trace}\,\big(\Theta^{\wedge k}\big)
		\end{equation}
		is a closed form on $Z$. Furthermore,  the cohomology class of ${\rm ch}(\nabla)$ (which is independent of the connection)   is the Chern character of $E$.
 
Let $\nabla'$ be a connection on $F'$. Then the Chern character of $i_!{\mathcal F}$ is the cohomology class defined by the closed form
\begin{equation}\label{chernnabla'}
{\rm ch}(\nabla'\otimes\Hat\nabla) -{\rm ch}(\nabla').
\end{equation}
The curvature of the tensor product of the connections is 
$$\Theta'\otimes I+I\otimes \Hat\Theta.$$
 On $X\setminus U$ this form reduces to $\Theta'\otimes I$. So, from (\ref{chernnabla'}), it follows that    the Chern character of $i_!{\mathcal F} $ can be  represented by differential  form that vanishes  out of $U$. The fact that this property holds for any neighbourhood of $Y$ indicates that ${\rm ch}(i_!{\mathcal F})$
 is localized on $Y$.
	
	\smallskip
{\it (b) Localization via gauge fields.}
Let $Z$ be a compact projective variety  and consider an embedding $j:Z\hookrightarrow M$ of $Z$ in a smooth variety $M$.  So, the Poincar\'e duality allows us to identify   homology and cohomology of $M$. Furthermore,   if ${\mathcal F}$ is a coherent sheaf over $Z$, then $j_!{\mathcal F}$ admits a resolution $E_{\bullet}$ by  locally free sheaves. The Chern character ${\rm ch}(E_{\bullet})\in H_*M$ restricts to zero in $H_*(M\setminus Z)$. Thus,  ${\rm ch}(E_{\bullet})$ is the image  of an element of $H_*Z$ which is denoted  by ${\rm ch}^M_Z({\mathcal F})$.
		
		Next, we show a construction of the localized class ${\rm ch}^M_Z({\mathcal F})$ by means of  strength of gauge fields over $M$, 
		developing what has already been sketched in the Introduction. 
		
		Let 
		\begin{equation}\label{resolutionjF}
		 0\to E_r\to E_{r-1}\to\dots\to E_0\to j_!{\mathcal F} \to 0
		 \end{equation}
		be a locally free resolution of $j_!{\mathcal F}$.
		We denote by $E_{\bullet}$ the complex 
		$0\to E_r\to E_{r-1}\to\dots\to E_0\to 0.$
			As $j_!{\mathcal F}$ is zero out $Z$, one has the following exact sequence
		$$(E_{\bullet}\big|_{M\setminus Z},\,d_{\bullet}):\; 0\to E_r\big|_{M\setminus Z}\overset{d_r}{\longrightarrow} E_{r-1}\big|_{M\setminus Z}\overset{d_{r-1}}{\longrightarrow}\dots\to E_0\big|_{M\setminus Z}\to 0.$$
	
It is possible to construct a connection   $\tilde\nabla_i$ on the vector bundle $ E_i\big|_{M\setminus Z}$, such that the  family $\{\tilde\nabla_i\}_{i=0}^r$ is compatible with the  operators $d_i$ (see \cite[Lemma 4.17]{B-B}).  We set
$${\rm ch}(\tilde\nabla_{\bullet}):=\sum_{i=0}^r(-1)^i{\rm ch}(\tilde\nabla_i),$$
 where ${\rm ch}(\tilde\nabla_i)$ is defined according to (\ref{chern-Nabla}).  The mentioned compatibility implies ${\rm ch}(\tilde\nabla_{\bullet})=0$ (see \cite[Lemma 4.22]{B-B}). 

Let $\nabla_{\bullet}$ be a family of connections for the complex $E_{\bullet}$, then ${\rm ch}(\nabla_{\bullet})$ is a closed form on $M$, which defines the Chern character of $j_!{\mathcal F}$,  by the exactness of (\ref{resolutionjF}). 

Since ${\rm ch}(\nabla_{\bullet}|_{M\setminus Z} )$ and ${\rm ch}(\tilde\nabla_{\bullet})$  define the same cohomology class on $M\setminus Z$, it turns out that ${\rm ch}(\nabla_{\bullet}|_{M\setminus Z} )=d\lambda$, with $\lambda$ a differential form on $M\setminus Z$. One may take  as $\lambda$  the Bott difference form constructed from $\nabla_{\bullet}|_{M\setminus Z}$ and $\tilde\nabla_{\bullet}$.

In summary, ${\rm ch}(\nabla_{\bullet})$ is a closed form on $M$ and its restriction to $M\setminus Z$ is equal to the differential of $\lambda$. Thus,  the pair $({\rm ch}(\nabla_{\bullet}),\, \lambda)$ defines an element of the relative cohomology $H^*(M,\,M\setminus Z)$ \cite[page 78]{B-T}. By  the Alexander-Lefschetz duality isomorphism $H^*(M,\,M\setminus Z)\simeq H_*(Z)$, that pair determines an element in $H_*(Z)$
which is the localization   of ${\rm ch}(j_!{\mathcal F})$.

\smallskip	
		
{\it (c) Construction of Baum-Fulton-MacPherson.}
Next, we summarize the construction of  ${\rm ch}^M_Z({\mathcal F})$ 
given in \cite{B-F-M}. Let $(E_{\bullet}, \,d_{\bullet})$ be the above complex  of locally free sheaves on $M$ obtained from a resolution of $j_!{\mathcal F}$. So, this complex is exact   on $M\setminus Z$. Then on $M\setminus Z$ we can fix splittings for each $i$ 
$${\rm ker}(d_i)|_{M\setminus Z}\oplus {\rm ker}(d_{i-1})|_{M\setminus Z}\simeq E_i|_{M\setminus Z}.$$
Hence, over $M\setminus Z$ one has the following isomorphism 
$$\beta:E_{\rm even}=\bigoplus_i  E_{2i}\overset{\sim}{\longrightarrow} \bigoplus_i\,{\rm ker}(d_i) \overset{\sim}{\longrightarrow}\bigoplus_i E_{2i+1} =E_{\rm odd}.$$
The difference bundle $d(E_{\bullet}):=d(E_{\rm even},\,E_{\rm odd},\,\beta)$ \cite{A-H} is an element of  the relative Grothendieck group  $ K^0(M,\,M\setminus Z)$, and its Chern character determines by duality an element ${\rm ch}^M_Z(E_{\bullet})\in H_*(Z,\,{\mathbb Q})$.
	
	
	{\it Definition \cite{B-F-M}.}	
Given a coherent sheaf ${\mathcal F}$ over $Z$, an embedding $j:Z\to M$, where $M$ is a smooth variety, and   a locally free resolution $E_{\bullet}$  of $j_!{\mathcal F}$. One defines 
$$\tau({\mathcal F}):={\rm td}(M)\cap {\rm ch}^M_Z(E_{\bullet}),$$  
where  $\cap$ denotes the cap product between cohomology and homology of $M$.

\smallskip

  In   \cite{B-F-M} it is proved that $\tau({\mathcal F})$ 
	 is independent of the  embedding   $j$ and of the resolution  $E_{\bullet}$. 
	 Furthermore,   $\tau$ determines a natural transformation between the functors $K_0$ and $H_*$ (the singular homology with rational coefficients) $\tau:K_0\to H_*$ such that
	\begin{equation}\label{GRR}
	\tau\big( [E] \cdot [{\mathcal F}]  \big)=({\rm ch} [E] )\cap(\tau[{\mathcal F}]),
	\end{equation}
		 for any $[E]\in K^0Z$, $[{\mathcal F}]\in K_0Z$.
		
		The  natural transformation $\tau$ 
		when $Z$ is a {\em smooth} variety, reduces to 
\begin{equation}\label{smooth}
[{\mathcal H}]\in K_0(Z)\mapsto P({\rm ch}({\mathcal H}) \,{\rm td}(Z))\in H_*(Z,\,{\mathbb Q}),
\end{equation} 
where $P$ denotes the Poincar\'e isomorphism $H^*Z\to H_*Z$.
		
			Let $Y$ be  a  local complete intersection    in a smooth variety $X$, we denote by $N_{Y/X}$ the normal bundle to $Y$ in $X$.
		The virtual vector bundle $T_Y$ of $Y$ is the element $T_X|_Y-N_{Y/X}$ of $K^0(Y)$, where $T_X$ is the tangent bundle of $X$. In this simple case, 
		\begin{equation}\label{taumathOY}
		\tau({\mathcal O}_Y)={\rm td}(T_Y)\cap[Y], 
		\end{equation}
		where ${\rm td}(T_Y)$ is the Todd class of the virtual bundle $T_Y$   (see \cite[Chapter IV]{B-F-M}).

	\smallskip	

		\subsection{Definition of $Q$.}
		
		Given a brane ${\mathcal G}\in D^b(Z)$ on the variety $Z$, we propose to define its charge $Q({\mathcal G})$ as the value of $\tau$ on the image of ${\mathcal G}$ by the map (\ref{DmapsstoK}); that is,
		\begin{equation}\label{Q}
		Q({\mathcal G}):=\tau([{\mathcal G}]).
		 \end{equation}

		Let $f:Z\to X$ be a proper morphism between varieties, then one has the following commutative diagram
	$$	\xymatrix{ D^b(Z)\ar[d]^{Rf_*} \ar[r]^{\alpha_Z} &  K_0Z\ar[d]^{f_*} \ar[r]^{\tau} & H_*Z \ar[d]^{f_*} \\
		D^b(X)\ar[r]^{\alpha_X} & K_0X \ar[r]^{\tau} & H_*X\, ,  }\,$$
		where $f_*$ on the right  is the map induced by $f$ on the homologies, and the horizontal arrows of the left square are defined in (\ref{DmapsstoK}).
		This diagram is commutative, since $\tau$ is a natural transformation between the functors $K_0$ and $H_*$ and by the definition of the map (\ref{DmapsstoK}).  
	
			\begin{Prop}\label{P:Coherencia}
Let $i:Y\hookrightarrow X$ be a closed subvariety of $X$ and ${\mathcal F}$ be a coherent sheaf on $Y$, then
$$ i_* Q(\mathcal {\mathcal F} )=Q(i_!{\mathcal F}):$$
\end{Prop}
{\it Proof.} Let us consider the above commutative diagram when 
		  $Z=Y$   and $f=i$. Since $i$ is a closed embedding, $Ri_*{\mathcal F} =i_!{\mathcal F}$. From   the equality $ i_* \circ\tau\circ\alpha_Y=\tau\circ\alpha_X\circ Ri_*$ applied to the element
			${\mathcal F}\in D^b(Y)$ together with (\ref{Q}), it follows 
			the proposition. \qed
		
		\smallskip	
	Let $Y$ be a local complete intersection of a nonsingular variety $X$ defined by a section of the algebraic vector bundle $N$ on $X$, and
let ${\mathcal E}$ be an object of $D^b(Y)$ consisting of locally free sheaves. By definition of $Q$ together with  the fact that $[{\mathcal O}_Y]$ is the identity of $K^0Y$
$$Q({\mathcal E})=\tau(\alpha_Y({\mathcal E}))=\tau(\alpha_Y({\mathcal E})\cdot [{\mathcal O}_Y]).$$
Since ${\mathcal E}$ is a complex of vector bundles,   from (\ref{GRR})  together with (\ref{taumathOY}) one deduces
 $$\tau(\alpha_Y({\mathcal E})\cdot [{\mathcal O}_Y])={\rm ch}({\mathcal E})\cap\big({\rm td}[T_Y]\cap [Y]\big).$$
From the equality in $K^0Y$, $[T_Y]=[i^*(T_X)]-[i^*N]$, it follows
\begin{equation}\label{Qmathcal E}
Q({\mathcal E})= \big({\rm ch}({\mathcal E})\,i^*({\rm td}(X)\, {\rm td}(N)^{-1}) \big)\cap [Y],
\end{equation}
		where $i$ is the inclusion $Y\hookrightarrow X$.
		
		\begin{Prop}\label{Prop:Charge}
		 If $Y$ is a local complete intersection of a nonsingular variety $X$ defined by a section of the algebraic vector bundle $N$ on $X$, and
 ${\mathcal E}$   an object of $D^b(Y)$ consisting of locally free sheaves, then   the charge of the brane ${\mathcal E}$ is given by (\ref{Qmathcal E}).
		\end{Prop}


\section{Charge of branes on a toric variety}\label{S:Euler}
\noindent
 Let $\Delta$ be a lattice polytope of dimension $n$ in $({\mathbb R}^n)^*$. Given a facet $A$ of $\Delta$, there is   a unique vector $v_A\in{\mathbb Z}^n$ conornal to $A$ and inward to $\Delta$. So, $A$ is on  the hyperplane in $({\mathbb R}^n)^*$ of equation
 $\langle m,\,v_A\rangle=-\kappa_A$, with $\kappa_A\in {\mathbb Z}$ and
$$\Delta=\bigcap_{A\, {\rm facet}}\{  m\in ({\mathbb R}^n)^*\,|\, \langle m,\,v_A\rangle\geq -\kappa_A   \}.$$

We denote by $X$ the toric variety determined by $\Delta$ and by $D_A$ the divisor of $X$ associated to the facet $A$. The divisor 
$$D_{\Delta}:=\sum_{A\, {\rm facet}} \kappa_AD_A$$
 is Cartier, ample and base pointfree \cite[page 269]{C-L-S}. Moreover, the torus invariant divisor $K_X:=-\sum_A D_A$ is a canonical divisor of $X$.   

From now on, we assume that $\Delta$ is a {\it reflexive} polytope; that is,   $\kappa_A=1$ for all $A$. 
Then the canonical sheaf $\omega_X={\mathcal O}_X(K_X)$ is a line bundle and $X$ is Fano (that is, the anticanonical divisor $-K_X$ is ample).

On the other hand, a generic element $Y$ of the linear system $|-K_X|$ is a Calabi-Yau orbifold 
 (for  details see \cite[Sec. 4.1]{C-K}).
 Henceforth, $Y$ will be a generic anticanonical hypersurface of $X$ and we denote by $i$ the inclusion $Y\hookrightarrow X$. 

 We set $N$ for the line bundle on $X$ associated with ${\mathcal O}_X(-K_X)$. Then $Y$ is a local complete intersection defined by a section of $N$. Moreover, the first Chern class of $X$ is also equal to $c_1(N)$ (\cite[page 625]{C-L-S}). We set 
 \begin{equation}\label{def:a}
a:={ c}_1(X)=c_1(N).
\end{equation}

The multiple $-(n-1)K_X$ of the anticanonical divisor $-K_X$ is a very ample divisor on $X$ \cite[page 71]{C-L-S}. We set ${L}$ for the restriction to $Y$ of the very ample invertible sheaf on $X$ defined by the Cartier divisor $-(n-1)K_X$. In \cite{Vina18}, we showed that   ${L}$   generates a spanning class in the bounded derivative category $D^b(Y)$.


\smallskip

{\it Arithmetic genus.}
Now let us assume that $X$ 
is smooth and ${\rm dim}_{\mathbb C}X=4$. So, $Y$ is a Calabi-Yau orbifold of dimension $3$. Our purpose is to express the arithmetic genus $\chi(Y,\, E)$, where $E$ is the restriction to $Y$ of a line bundle $E'$ over $X$,
in terms of  intersections of divisors in $X$.

We will apply the generalization of the Hirzebruch-Riemann-Roch formula to singular varieties  \cite{B-F-M}.
In our case,
\begin{equation}\label{HRR}
\chi(Y,\,{E})=\epsilon\big( {\rm ch}(E)\cap\tau ({\mathcal O}_Y)  \big),
\end{equation}
where $\epsilon:H_*Y\to{\mathbb Q}$ is the map induced on the rational homology of $Y$ by the mapping $Y\to\{\rm pt\}$, and $\tau({\mathcal O}_Y)$ is given  by (\ref{taumathOY}).
 
The Chern class of the virtual tangent bundle  $T_Y$ is
$${\rm c}(T_Y)=i^*\big({\rm c}(T_X)({\rm c}(N))^{-1}\big)\in H^*Y.$$
 Taking into account the notation $a={\rm c}_1({\mathcal O}_X(-K_X)),$ 
$${\rm c}(N)^{-1}=\sum_{k\geq 0}(-1)^ka^k.$$

By (\ref{def:a}), the Chern class of the tangent bundle $T_X$ can be written
$${\rm c}(T_X)=1+a+{\rm c}_2(X)+{\rm c}_3(X)+{\rm c}_4(X).$$
Hence,
$${\rm c}(T_X){\rm c}(N)^{-1}=1+ {\rm c}_2(X)+({\rm c}_3(X)-a{\rm c}_2(X))+\dots.$$
That is,
$${\rm c}_1(T_Y)=0,\;\; {\rm c}_2(T_Y)=i^*{\rm c}_2(X),\;\; {\rm c}_3(T_Y)=i^*({\rm c}_3(X)-a{\rm c}_2(X)).$$ 

For the Todd class ${\rm td}(T_Y)$, we have
\begin{align}\notag
{\rm td}(T_Y) &=1+\tfrac{1}{2}{\rm c}_1(T_Y)+\tfrac{1}{12}({\rm c}_2(T_Y)+{\rm c}_1^2(T_Y)) + \tfrac{1}{24}{\rm c}_1(T_Y){\rm c}_2(T_Y) \\
&=i^*\big( 1+ \tfrac{1}{12}c_2(X) \big) \notag 
\end{align}

The Chern character of $E=i^*E'$ is
 $${\rm ch}(E)=i^*\sum_{k\geq 0}\frac{(c_1(E'))^k}{k!}.$$

The summands of degree six in ${\rm ch}(E)\, {\rm td}(T_Y)\in H^*(Y)$ are  
$$i^*\big(\tfrac{1}{6}(c_1(E')^3)+\tfrac{1}{12}(c_1(E') c_2(X)\big).$$
 By (\ref{HRR}),
\begin{equation}\label{chi(Y,E)}
\chi(Y,\,E)=i^*\big(\tfrac{1}{6}(c_1(E')^3)+\tfrac{1}{12}(c_1(E') c_2(X)\big)\cap [Y].
 \end{equation}

Given a face $B$ of $\Delta$, we set $V(B)$ for the closure of the orbit in $X$ determined by $B$ through the toric action. $[V(B)]$ will denote the corresponding homology class. Some of those classes are involved in the formula  for the arithmetic genus $\chi(Y,\, E)$, according to  the following theorem.

 \begin{Thm}\label{Thm:3} Under the above hypotheses, if the line bundle $E'$ is determined by the divisor $D$, then
 $$\chi(Y,\,E)=-[K]\cdot\Big(\tfrac{1}{6}[D]^3 + \tfrac{1}{12}[D]\cdot  \sum_{B\in{\sf F}_2} [V(B)] \Big),$$ 
 where ${\sf F}_2$ is the set of codimension $2$ faces of $\Delta$, $K$ is the canonical divisor of $X$, and  $[\ast]$ (for $\ast=K, D$) denotes the homology class of $X$ defined by the divisor $\ast$. 
\end{Thm}
 {\it Proof.}
Since second Chern class of $X$
is  \cite[page 625]{C-L-S}
$$c_2(X)=P^{-1}\big(\sum_{B\in{\sf F}_2} [V(B)] \big),$$
 the theorem follows from (\ref{chi(Y,E)}).
 \qed

\smallskip

In the particular case that $E$ is the line bundle $L$ mentioned above, then one can take $D=-3K_X$ and we have the following corollary.
 \begin{Cor}\label{Cor:L}
 The arithmetic genus of $L$ is
 $$\chi(Y,\,L)=[K]^2\cdot\Big(\tfrac{9}{2}[K]^2 + \tfrac{1}{4} \sum_{B\in{\sf F}_2} [V(B)] \Big).$$
\end{Cor}


\smallskip
{\it Charge of the brane defined by a divisor.}
Let $F'$ be the line bundle on $X$ defined by a divisor $D$ of $X$, and let 
us  we assume that $X$ is  smooth of dimension $3$. 
As in Section
\ref{Ss:Consistence}, $F'$ determines a coherent sheaf ${\mathcal F}$ on $Y$. 
We will calculate the image of $Q({\mathcal F})$ in $H_*X$; that is, $i_*Q({
\mathcal F})$. According to the consistency condition (\ref{i_*Q}), we can 
calculate $Q(i_!{\mathcal F})$. Since $X$ is smooth, by (\ref{smooth})
 $$Q(i_!{\mathcal F})=P\big({\rm td}(X)\,{\rm ch}(i_!{\mathcal F})  \big).$$
On the other hand, denoting $a:=c_1({\mathcal O}(-K_X))$ as above,
$${\rm td}(X)=1+\tfrac{1}{2}a+\tfrac{1}{12}(a^2+{\rm c}_2(X)) +\dots.$$
 By (\ref{rmch})
$${\rm ch}(i_!{\mathcal F})=a\big(1+\tfrac{1}{2} a+\tfrac{1}{6} a^2 \big) \big(1+{\rm c}_1(F')+\tfrac{1}{2}{\rm c}_1(F')^2\big).$$
Thus,
\begin{align}\notag
&{\rm td}(X)\,{\rm ch}(i_!{\mathcal F})=\\ \notag
 &  a\Big(1+a+{\rm c}_1(F')+\tfrac{1}{2}\big(a^2+{\rm c}_1(F')^2+\tfrac{1}{6}{\rm c}_2(X)\big)  +a{\rm c}_1(F')     \Big)\in H^*X.\notag
\end{align}
Hence, the charge of ${\mathcal F}$ considered as an element of $H_*X$ can be written in terms of $D$ and the canonical divisor $K$ of $X$.
\begin{equation}\label{tocho}
-[K]\cdot\Big( 1-[K] +[D] +\tfrac{1}{2}\big([K]^2 +[D]^2+ \tfrac{1}{6} \sum_{B\in{\sf F}_2} [V(B)]\big) -[K]\cdot [D]                         \Big),
\end{equation}
where   the notation above introduced is used.

 \begin{Prop}\label{Prop:charge}
 Let ${\mathcal F}$  be the sheaf on $Y$ defined by restriction of the line bundle $F'$ on  $X$. If ${\rm dim}_{\mathbb C}  X=3$ and  $F'$ is determined by the divisor $D$, then the homology class $i_*Q({\mathcal F})$ is given by (\ref{tocho}).
\end{Prop}

The divisor $-2K_X$ determines the above line bundle  ${L}$, generator of a spanning class  of $D^b(Y)$. Thus, we have the following corollary.
\begin{Cor}\label{Cor:mathacalL}
Let ${\mathcal L}$ be the brane on $Y$ defined by the divisor $-2K$. Then
$$i_*Q({\mathcal L})=-[K]\cdot\Big(1-3[K]+\tfrac{9}{2}[K]^2+\tfrac{1}{12} \sum_{B\in{\sf F}_2} [V(B)]   \Big).$$
\end{Cor}

In the case that  $X$ is a complex surface,  in the expression for $i_*Q({\mathcal F})$ reduces to the first three terms of (\ref{tocho}).

\begin{Cor}\label{Coro:charge} Let ${\mathcal F}$  be the sheaf on $Y$ defined by restriction of the line bundle $F'$ on the surface $X$. If $F'$ is determined by the divisor $D$, then the homology class $i_*Q({\mathcal F})$ is given by 
$$i_*Q({\mathcal F})=-[K]+[K]\cdot([K]-[D]).$$
\end{Cor}


\end{document}